\documentclass[11pt]{amsart}

\usepackage{amsmath}
\usepackage{amsthm}
\usepackage{amsfonts}
\usepackage{amssymb,latexsym}
\usepackage[all]{xy}
\usepackage{graphicx}
\usepackage[mathscr]{eucal}
\usepackage{verbatim}
\usepackage{hyperref}

\theoremstyle{plain}
    \newtheorem{thm}{Theorem}[section]

    \newtheorem{cor}[thm]{Corollary}

\theoremstyle{definition}

\theoremstyle{remark}

\numberwithin{equation}{section}

\newcommand{\fq}{\mathbb{F}_{q}}

\begin{document}

\title[Counting irreducible polynomials over finite fields]{Counting irreducible polynomials over finite fields using the inclusion-exclusion principle}
\date{\today}

\author{Sunil K. Chebolu}
\address{Department of Mathematics \\
Illinois State University \\
Normal, IL 61790, USA} \email{schebol@ilstu.edu}

\author{J\'{a}n Min\'{a}\v{c}}
\address{Department of Mathematics\\
University of Western Ontario\\
London, ON N6A 5B7, Canada}
\email{minac@uwo.ca}


\begin{abstract}
 C. F. Gauss discovered a beautiful formula for the number of irreducible polynomials of a given degree over a finite field.  Using just very basic knowledge of finite fields and the inclusion-exclusion formula, we show how one can see the shape of this formula and its proof almost instantly.
\end{abstract}

\maketitle
\thispagestyle{empty}

Why  there are exactly
\[\frac{1}{30} (2^{30} - 2^{15} - 2^{10} - 2^{6} + 2^{5} + 2^{3} + 2^{2} -2)\]
irreducible monic polynomials of degree 30 over the field of two elements? In this note we will show how one can see the answer instantly using just very basic knowledge of finite fields and the well-known  inclusion-exclusion principle.

To set the stage, let $\fq$ denote the finite field of $q$ elements.
Then in general, the  number of monic irreducible polynomials of degree $n$ over the finite field $\fq$ is given by
Gauss's formula
\[  \frac{1}{n} \sum_{d | n} \mu(n/d) q^d, \]
where  $d$ runs over the set of all positive divisors of $n$ including $1$ and $n$, and  $\mu(r)$ is the M\"obius function. (Recall that $\mu(1) =1$ and $\mu(r)$ evaluated at a product of distinct primes is 1 or -1 according to whether the number of factors is even or odd. For all other natural numbers $\mu(r) =0$.)
This beautiful formula is well-known and was discovered by Gauss \cite[p. 602-629]{Gauss} in the case when $q$ is a prime. 

 We will present a proof of this formula that uses only elementary facts about finite fields and the inclusion-exclusion principle. Our approach offers the reader a new insight into this formula because our proof gives a precise field theoretic meaning to each summand in the above formula.  The classical proof  \cite[p. 84]{IR} which uses the M\"obius' inversion formula does not offer this insight. Therefore we hope that students and users of finite fields may find our approach helpful.
It is surprising that our simple argument  is not available in textbooks, although it must be known to some specialists.

\vskip 6mm
\noindent
\textbf{Proof of Gauss's formula} \ \  Before we present our proof we  collect some basic facts about finite fields that we will need. These facts and their proofs can be found in almost any standard algebra textbook that covers finite fields. See for example \cite[Chapter 14.3]{DummitFoote},
\cite[Chapter 7.1, 7.2]{IR}, or \cite[Chapter 20.1]{Judson}.

\begin{enumerate}
\item A finite field of order $q$ exists if and only if $q$ is a prime power. Moreover, such a field is unique up to isomorphism, and is denoted by $\fq$.
\item $\mathbb{F}_{q^n}$ is the splitting field of any irreducible polynomial $p(x)$ of degree $n$ over $\fq$. (This means $p(x)$ factors into linear factors over $\mathbb{F}_{q^n}$ but not over any smaller subfield of $\mathbb{F}_{q^n}$.)
\item The roots of an irreducible polynomial over $\fq$ are always distinct.
\item No two distinct irreducible polynomials over  $\fq$ can have a common root.
\item $\mathbb{F}_{q^a} \subseteq \mathbb{F}_{q^b}$ if and only $a$ divides $b$.
\end{enumerate}

With these basic facts under our belt we  proceed to show that the total number of irreducible monic polynomials of degree $n$ over $\fq$ is equal to
\[ \frac{1}{n} \sum_{d | n} \mu(n/d) q^d.\]
The case $n =1$ is easy because every degree one monic polynomial is irreducible.  In fact, the total number of degree one monic polynomials  over $\fq$ is equal to $q$, and this is
exactly what we get from the above formula when we plug in $n =1$.
Therefore for the rest of the proof we will assume that  $n > 1$.  Let $\mathcal{P}_n$ denote the collection of all irreducible monic
polynomials of degree $n$ over $\fq$ and  let $\mathcal{R}_n$ be the union of all the roots of all the polynomials in $\mathcal{P}_{n}$.  Note that fact (2) ensures that the roots thus obtained are contained in $\mathbb{F}_{q^n}$.
Moreover, using facts (3) and (4) we conclude that $\mathcal{R}_n$ is the disjoint union of $n$-element sets, one for each polynomial in $\mathcal{P}_n$.
Thus, 
\[|\mathcal{R}_n| = n|\mathcal{P}_n|. \]

Therefore, it is enough  to compute  $|\mathcal{R}_n|$. To this end, observe that
\begin{eqnarray*}
\mathcal{R}_n& = & \{ \alpha \text{ in } \mathbb{F}_{q^{n}} \,|\ \ [\mathbb{F}_{q}(\alpha) \colon \fq] = n \},\\
& = & \{ \alpha \text{ in } \mathbb{F}_{q^n} \, | \,  \alpha \text{ is not contained in any proper subfield of } \mathbb{F}_{q^n} \},\\
& = & \{ \alpha \text{ in } \mathbb{F}_{q^n} \, | \,  \alpha \text{ is not contained in any maximal proper subfield of
 } \mathbb{F}_{q^n} \}
 \end{eqnarray*}

Let $n = u^{a} v^{b} w^{c}  \cdots $ be the prime factorization of $n$ with $r$ distinct prime factors (recall that $n > 1$). Then the maximal subfields of $\mathbb{F}_{q^n}$ by fact (4) are of the form
\[F_u = \mathbb{F}_{q^{\frac{n}{u}}},  F_v = \mathbb{F}_{q^{\frac{n}{v}}},  F_w = \mathbb{F}_{q^{\frac{n}{w}}}, \cdots .\]
Then by the third interpretation of $\mathcal{R}_{n}$ given above, we have
\[ | \mathcal{R}_n | = | (F_u \cup F_v \cup F_{w} \cdots  )^c| \]
where the complement is taken in the field $\mathbb{F}_{q^n}$.
Fact (4) implies that the lattice of subfields of $\mathbb{F}_{q^n}$ that contain $\fq$ is isomorphic to the lattice of divisors of $n$.
Thus, we have   $F_u \cap F_v = \mathbb{F}_{q^{\frac{n}{u v}}}$ and $F_u \cap F_v \cap F_w = \mathbb{F}_{q^{\frac{n}{ u v w}}}$, etc.

The cardinality of $\mathcal{R}_{n}$ can now be computed using the inclusion-exclusion principle as follows.
\begin{align*}
|\mathcal{R}_n| =  \, & q^n \\
 &  -q^{\frac{n}{u}} - q^{\frac{n}{v}} -  q^{\frac{n}{w}}  - \cdots  \\
 &  + q^{\frac{n}{u v}} + q^{\frac{n}{u w}} + q^{\frac{n}{v w}}+  \cdots  \\
 & \ \ \ \ \ \  \ \ \ \ \cdots \\
 & + (-1)^r q^{\frac{n}{u v w \cdots }}.
\end{align*}

Finally, the formula for  $|\mathcal{P}_n|$ takes the desired form when we divide $|\mathcal{R}_{n}|$ by $n$ and
use the M\"obius function.

We end by pointing out that if one is interested in counting the cardinality of all irreducible polynomials of degree $n$ (not necessarily monic) over $\fq$ then we simply multiply $|\mathcal{P}_n|$ by $q-1$. This is because every such polynomial $p(x)$ can be uniquely written as $\alpha r(x)$ where
$\alpha$ is a non-zero element of $\fq$ and $r(x)$ is an irreducible monic polynomial of the same degree.

\vskip 5mm
\noindent
\textbf{Acknowledgements.}  We would like to thank our students at Illinois State University and University of Western Ontario for their inspiration and insistence on penetrating mysteries of finite fields. 
We also would like to thank Fusun Akman, Micheal Dewar, and Alan Koch for helping us improve the exposition through their interesting and very encouraging comments on this paper. Finally, we would like to thank the anonymous referee for 
their suggestions to improve the exposition. 

\vskip 10mm



\end{document}